\documentclass{article}
\catcode`\@=11
\newtheorem{theorem}{Theorem}
\newtheorem{lemma}{Lemma}
\newtheorem{proposition}{Proposition}

\newtheorem{example}{Example}

\newtheorem{corollary}{Corollary}
\def\demo{\noindent{\bf Proof .-}}
\def\section{\@startsection {section}{1}{\z@}{-3.5ex plus -1ex
minus-.2ex}{2.3ex plus .2ex}{\normalsize\bf}}

\newcommand{\rad}{{\rm{Rad}}}

\begin{document}
\begin{center}
{\Large\bf \textsc{On binomial equations defining rational normal scrolls}}
\end{center}
\vskip.5truecm
\begin{center}
{Margherita Barile\footnote{Partially supported by the Italian Ministry of Education, University and Research.}\\ Dipartimento di Matematica, Universit\`{a} di Bari, Via E. Orabona 4,\\70125 Bari, Italy\\barile@dm.uniba.it}
\end{center}
\vskip1truecm
\noindent
{\bf Abstract} We show that a rational normal scroll can in general be set-theoretically defined by a proper subset of the 2-minors of the associated two-row matrix. This allows us to find a class of rational normal scrolls that are almost set-theoretic complete intersections. 

\section*{Introduction}
If $I$ is an ideal in a commutative ring with identity $R$, then we say that the elements $P_1,\dots, P_s\in R$ {\it generate I up to radical} if $\rad(P_1,\dots, P_s)=\rad(I)$. Finding  such elements $P_i$ is especially interesting when $R$ is a polynomial ring over an algebraically closed field: in that case, according to Hilbert Nullstellensatz, the above condition means that the algebraic variety $V(I)$ is defined by $P_1=\cdots=P_s=0$. It is interesting to investigate those cases where we can take $s$ smaller than the minimum number of generators of $I$. Such a study has been made, e.g., for several classes of determinantal ideals associated with matrices of indeterminates: those associated with a generic matrix in \cite{BS}, those associated with matrices containing symmetries in \cite{BV}, \cite{B}, \cite{Va}, and \cite{V}. In this paper we consider the ideals generated by the 2-minors of the matrices composed of catalecticant blocks with two rows. These are the defining ideals of rational normal scrolls (see \cite{EG} or \cite{H}). We show that in most cases they can be  generated up to radical by a proper subset of the generating minors. We thus give a new class of binomial ideals which can be generated, up to radical, by a proper subset of the generating binomials. This will allow us to partially extend the main result in \cite{BV}; in particular we will find ideals $I$ defining rational normal scrolls which can be generated, up to radical, by height$(I)$ or height$(I)$+1 elements.  Other classes of binomial ideals sharing similar properties are the toric ideals described in \cite{B1}, \cite{B2}, \cite{B3}, \cite{BL}, \cite{BMT1} and \cite{BMT2}. 

\section{The statement of the main result}
Let $R$ be a commutative ring with identity. Let $r$ be a positive integer and consider the two-row matrix with entries in $R$
\begin{equation}\label{2}A=\left(B_1\vert B_2\vert\dots\vert B_r\right),\end{equation}
\noindent
where, for all $i=1,2,\dots, r$, $B_i$ is the $2\times c_i$ catalecticant matrix
$$B_i=\left(\begin{array}{ccccc}
x_{i0}&x_{i1}&\dots&x_{i{c_i-2}}&x_{i{c_i-1}}\\
x_{i1}&x_{i2}&\dots&x_{i{c_i-1}}&x_{i{c_i}}
\end{array}
\right).
$$
\noindent
 We will call $B_i$ the $i$th block of $A$. Let $I$ be the ideal of $R$ generated by the 2-minors of $A$. 
For all $i=1,\dots, r$, let $M_i$ denote the set of 2-minors of $B_i$; if $B_i$ consists of only one column, then we set $M_i=\emptyset$. Moreover, for all indices $i,j$ such that $1\leq i<j\leq r$, we set
$$N_{ij}=\{x_{i0}x_{jc_j}-x_{i1}x_{jc_j-1},\ x_{ic_i-1}x_{j1}-x_{ic_i}x_{j0}\};$$
\noindent
in other words, the elements of $N_{ij}$ are the 2-minors formed by the  first column of $B_i$ and the last column of $B_j$, and by the last column of $B_i$ and the first column of $B_j$ respectively. Our main goal is to prove the following result.
\begin{theorem}\label{main} Let $J$ be the ideal of $R$ generated by 
$$\big(\bigcup_{i=1}^rM_i\big)\bigcup\big(\bigcup_{1\leq i<j\leq r}N_{ij}\big).$$
\noindent Then 
$$\rad(I)=\rad(J).$$
\end{theorem} 
Note that the generating set of $J$ is a proper subset of the set of 2-minors of $A$ whenever at least two blocks of $A$ have more than one column. \newline
We remark that it suffices to prove Theorem \ref{main} for $r=2$. This will be done, in several steps, in the next section. Henceforth we will focus on the case $r=2$, which will simplify our argumentation and notation in a relevant way. In fact it is clear that it suffices to prove Theorem \ref{main} in the case where matrix $A$ has exactly two blocks. 

\section{The proof of the main result}
Consider the following matrix with entries in $R$: 
$$A_{c,d}=\left(\begin{array}{cccccccccc}
x_0&x_1&\dots&x_{c-1}&\vline&y_0&y_1&\dots&y_{d-1}\\
x_1&x_2&\dots&x_c&\vline&y_1&y_2&\dots&y_d
\end{array}\right),
$$
where $c,d$ are positive integers.  Let $I_{c,d}$ be the ideal of $R$ generated the 2-minors of $A_{c,d}$. Moreover,  for all $i$, $1\leq i\leq c$, let $S_i$ denote the set of all 2-minors of the submatrix of $A$ formed by the first $i$ columns; 
for all $j$, $1\leq j\leq d$, let $S_j$ denote the set of all 2-minors of the submatrix of $A$ formed by the first $j$ columns of the second block; finally,  let $J_{i,j}$ be the ideal of $R$ generated by
$$S_i\cup T_j\cup\{x_0y_j-x_1y_{j-1},\ x_{i-1}y_1-x_iy_0\}.$$

\noindent With respect to the notation we have just introduced, the claim of Theorem 1 in the case $r=2$ is 
$$\rad(I_{c,d})=\rad(J_{c,d}).$$
\noindent
Evidently $J_{c,d}\subset I_{c,d}$. Hence Theorem \ref{main} is proven once we have shown that 
\begin{equation}\label{claim}I_{c,d}\subset\rad(J_{c,d}).\end{equation}
\noindent
To this end we need two preliminary lemmas.
\begin{lemma}\label{lemma1} For all $i$, $2\leq i\leq c$,
$$x_1^{i-1}(x_0y_1-x_1y_0)^{i-1}\in J_{i,d}.$$
\end{lemma}
\demo We proceed by induction on $i$. It holds:
$$x_1(x_0y_1-x_1y_0)=x_0(x_1y_1-x_2y_0)+y_0(x_0x_2-x_1^2)\in J_{2,d},$$
\noindent which proves the claim for $i=2$. So let $i\geq3$ and suppose the claim true for $i-1$. We then have
\begin{eqnarray}\label{1} x_1(x_0y_1-x_1y_0)(x_{i-2}y_1-x_{i-1}y_0)&=&(x_{i-2}y_0y_1-x_{i-1}y_0^2)(x_0x_2-x_1^2)\nonumber\\
&-&x_0y_1^2(x_0x_{i-1}-x_1x_{i-2})\nonumber\\
&+&x_0y_0y_1(x_0x_i-x_2x_{i-2})\nonumber\\
&-&x_0y_0^2(x_1x_i-x_2x_{i-1})\nonumber\\
&+&x_0(x_0y_1-x_1y_0)(x_{i-1}y_1-x_iy_0)\nonumber\\
&\in& J_{i,d}.
\end{eqnarray}
\noindent
Since 
\begin{eqnarray*} J_{i-1, d}&=&(S_{i-1})+(T_d)+(x_0y_d-x_1y_{d-1})\\
&+&(x_{i-2}y_1-x_{i-1}y_0),
\end{eqnarray*}
\noindent
and $S_{i-1}, T_d, (x_0y_d-x_1y_{d-1})\subset J_{i,d}$, (\ref{1}) implies that 
$$ x_1(x_0y_1-x_1y_0)J_{i-1,d}\subset J_{i,d}.$$
\noindent
Thus, by induction, 
\begin{eqnarray*} x_1^{i-1}(x_0y_1-x_1y_0)^{i-1}&=&x_1(x_0y_1-x_1y_0)x_1^{i-2}(x_0y_1-x_1y_0)^{i-2}\\
&\in&x_1(x_0y_1-x_1y_0)J_{i-1, d}\subset J_{i,d},
\end{eqnarray*}
\noindent as required. This completes the proof.
\par\medskip\noindent
\begin{lemma}\label{lemma2}
For all $j$, $2\leq j\leq d$, 
$$J_{c, j-1}\subset\rad(J_{c,j}).$$
\end{lemma}
\demo
We have that $S_c, T_{j-1}, x_{c-1}y_1-x_cy_0\in J_{c,j}$. It remains to show that 
\begin{equation}\label{ast} x_0y_{j-1}-x_1y_{j-2}\in\rad(J_{c,j}).
\end{equation}
\noindent First note that 
\begin{eqnarray}\label{3} y_{j-1}(x_0y_{j-1}-x_1y_{j-2})&=&-x_0(y_{j-2}y_j-y_{j-1}^2)\nonumber\\
&+&y_{j-2}(x_0y_j-x_1y_{j-1})\in J_{c,j}.
\end{eqnarray}
\noindent
Moreover, for $j\geq3$, one has that
\begin{eqnarray}\label{4}
y_{j-2}(x_0y_{j-1}-x_1y_{j-2})&=&-x_0(y_{j-3}y_j-y_{j-1}y_{j-2})\nonumber\\
&+&y_{j-3}(x_0y_j-x_1y_{j-1})\nonumber\\
&+&x_1(y_{j-3}y_{j-1}-y_{j-2}^2)\in J_{c,j},
\end{eqnarray}
\noindent
whereas, by Lemma \ref{lemma1}, for $j=2$ we have
\begin{equation}\label{5} x_1^{c-1}(x_0y_{j-1}-x_1y_{j-2})^{c-1}=x_1^{c-1}(x_0y_1-x_1y_0)^{c-1}\in J_{c,j}.
\end{equation}
\noindent
For $j=2$, (\ref{3}) and (\ref{5}) imply that $y_{j-1}(x_0y_{j-1}-x_1y_{j-2}), x_1(x_0y_{j-1}-x_1y_{j-2})\in\rad(J_{c,j})$. It follows that $(x_0y_{j-1}-x_1y_{j-2})^2\in \rad(J_{c,j})$, which is equivalent to (\ref{ast}).  For $j\ge 3$ claim (\ref{ast}) follows similarly from (\ref{3}) and (\ref{4}).  This completes the proof.
\par\medskip\noindent

We are now able to prove (\ref{claim}). We proceed by induction on $c+d$. If $c+d=2$, then $c=d=1$ and there is nothing to prove. So let $c+d>2$ and suppose that the claim is true for all smaller values of $c+d$. Up to exchanging the blocks of $A_{c,d}$ we may assume that $d\geq2$.   Set
$$U_i=x_iy_d-x_{i+1}y_{d-1}\qquad\mbox{for }0\leq i\leq c-1,$$
\noindent and 
$$V_i=y_iy_d-y_{i+1}y_{d-1}\qquad\mbox{for }0\leq i\leq d-2.$$
\noindent
Then 
$$ I_{c,d}=I_{c, d-1}+(U_0,\dots, U_{c-1})+(V_0, \dots, V_{d-2}).
$$
\noindent
By induction $I_{c,d-1}\subset\rad(J_{c, d-1}).$ Hence
\begin{eqnarray*} I_{c,d}&\subset& \rad(J_{c,d-1}+(U_0,\dots, U_{c-1})+(V_0,\dots, V_{d-2}))\\
&=& \rad((S_c)+(T_{d-1})+(x_0y_{d-1}-x_1y_{d-2})+(x_{c-1}y_1-x_cy_0)\\
&&+(U_0)+(U_1,\dots, U_{c-1})+(V_0,\dots, V_{d-2}))\\
&=&\rad((S_c)+(T_d)+(x_0y_{d}-x_1y_{d-1})+(x_{c-1}y_1-x_cy_0)\\
&&+(x_0y_{d-1}-x_1y_{d-2})+(U_1,\dots, U_{c-1}))\\
&=&\rad(J_{c,d}+(x_0y_{d-1}-x_1y_{d-2})+(U_1,\dots, U_{c-1}))\\
&=&\rad(J_{c,d}+(U_1,\dots, U_{c-1})),
\end{eqnarray*}
\noindent
where the last equality follows from (\ref{ast}). If $c=1$, we have the required claim (\ref{claim}). So assume that $c\geq 2$. All we have to show is that 
for all indices $i$ with $1\leq i\leq c-1$, 
\begin{equation}\label{6} U_i\in\rad(J_{c,d}).\end{equation}
\noindent
We prove (\ref{6}) by induction on $i$. We have that
\begin{equation}\label{a}x_1(x_1y_d-x_2y_{d-1})=-(x_0x_2-x_1^2)y_d+x_2(x_0y_d-x_1y_{d-1})\in J_{c,d}.\end{equation}
\noindent Moreover 
$$y_{d-1}(x_1y_d-x_2y_{d-1})=(y_{d-2}y_d-y_{d-1}^2)x_2+y_d(x_1y_{d-1}-x_2y_{d-2}),$$
\noindent
where the first summand belongs to $J_{c,d}$ by definition and the second summand belongs to $I_{c, d-1}$ and, consequently, to $\rad(J_{c,d-1})$ 
by induction on $c+d$. Therefore, in view of Lemma \ref{lemma2},   
\begin{equation}\label{b}y_{d-1}(x_1y_d-x_2y_{d-1})\subset\rad(J_{c,d}).\end{equation}
\noindent
By (\ref{a}) and (\ref{b}) we conclude that
$$x_1y_d-x_2y_{d-1}\in\rad(J_{c,d}),$$
\noindent which proves our claim (\ref{6}) for $i=1$. Now take $i>1$ and suppose that (\ref{6}) is fulfilled for $i-1$. Then 
\begin{eqnarray}\label{7}
x_i(x_iy_d-x_{i+1}y_{d-1})&=&-(x_{i-1}x_{i+1}-x_i^2)y_d+x_{i+1}(x_{i-1}y_d-x_iy_{d-1})\nonumber\\
&\in&\rad(J_{c,d}),
\end{eqnarray}
\noindent because $x_{i-1}x_{i+1}-x_i^2\in J_{c,d}$, $U_{i-1}=x_{i-1}y_d-x_iy_{d-1}\in\rad(J_{c,d})$ by induction on $i$, and
\begin{eqnarray}\label{8}y_{d-1}(x_iy_d-x_{i+1}y_{d-1})&=&(y_{d-2}y_d-y_{d-1}^2)x_{i+1}+y_d(x_iy_{d-1}-x_{i+1}y_{d-2})\nonumber\\
&\subset&\rad(J_{c,d}).
\end{eqnarray}
\noindent
The inclusion holds because the first summand belongs to $J_{c,d}$ by definition and the second summand belongs to $I_{c, d-1}$ and, consequently, to $\rad(J_{c,d-1})$ by induction on $c+d$, and hence to $\rad(J_{c,d})$ by virtue of Lemma \ref{lemma2}. 
Claim (\ref{6}) then follows from (\ref{7}) and (\ref{8}). This completes the proof of (\ref{claim}).
 
\section{A consequence in positive characteristic}
According to Theorem \ref{main}, whose proof we have just accomplished, for $c,d\geq 2$, ideal $I_{c,d}$ can be generated, up to radical, by ${c\choose 2}+{d\choose 2}+2$ minors of $A_{c,d}$. This number can be lowered by one, in certain cases, if $R$ is a domain of positive characteristic. This is what we are going to show next, after proving a preliminary result.  Let $J_d=(S_d)+(T_d)+(x_0y_d-x_dy_0)$.
\begin{lemma}\label{lemma3} For all $d\geq 1$,
\begin{eqnarray*} x_0^dy_d^d&\equiv& x_1^dy_{d-1}^d\quad({\rm mod}\ J_d);\\
x_d^dy_0^d&\equiv& x_{d-1}^dy_1^d\quad({\rm mod}\ J_d).
\end{eqnarray*}
\end{lemma}
\demo We first prove that, for all $k=0,\dots, d-1$, 
\begin{equation}\label{claim2} x_0^{d-1}x_d\equiv x_1^kx_{d-k}x_0^{d-k-1}\quad(\mbox{mod }(S_d)).\end{equation}
\noindent
The claim is trivial for $k=0$. Suppose that $k\geq1$, and that the claim is true for $k-1$. Then 
\begin{eqnarray*} x_1^kx_{d-k}x_0^{d-k-1}&=& x_1^{k-1}x_1x_{d-k}x_0^{d-k-1}\\
&\equiv&x_1^{k-1}x_0x_{d-k+1}x_0^{d-k-1}\quad(\mbox{mod }(S_d)),\\
&=&x_1^{k-1}x_{d-k+1}x_0^{d-k}
\end{eqnarray*}
\noindent
since $x_0x_{d-k+1}-x_1x_{d-k}\in S_d$. By induction 
$$x_1^{k-1}x_{d-k+1}x_0^{d-k}\equiv x_0^{d-1}x_d\quad(\mbox{mod }(S_d)),$$
\noindent
which completes the induction step and the proof of (\ref{claim2}). In particular, (\ref{claim2}) for $k=d-1$ gives
\begin{equation}\label{1a} x_0^{d-1}x_d\equiv x_1^d\quad(\mbox{mod }(S_d)).\end{equation}
\noindent
Similarly one can show that
\begin{equation}\label{2a} y_d^{d-1}y_0\equiv y_{d-1}^d\quad(\mbox{mod }(T_d)).\end{equation}
\noindent
Now \begin{eqnarray*} x_0^dy_d^d&\equiv& x_0^{d-1}y_d^{d-1}x_0y_d\\
&\equiv& x_0^{d-1}y_d^{d-1}x_dy_0\quad(\mbox{mod }J_d)\\
&=&x_0^{d-1}x_dy_d^{d-1}y_0\\
&\equiv&x_1^dy_{d-1}^d\quad(\mbox{mod }J_d),
\end{eqnarray*}
\noindent
where the last relation follows from (\ref{1a}) and (\ref{2a}).
This proves the first part of the claim. The second part follows by symmetry.
\par\medskip\noindent
\begin{corollary}\label{corollary4} Let $p$ be a prime, and suppose that $R$ is a domain of characteristic $p$. Then, for every positive integer $h$, 
$$\rad(I_{p^h, p^h})=\rad(J_{p^h}).$$
\end{corollary}
\demo On the one hand we have
\begin{equation}\label{4a} J_{p^h}\subset I_{p^h, p^h}, \end{equation}
\noindent
on the other hand, from  (\ref{claim}) it follows that
\begin{equation}\label{3a} I_{p^h, p^h}\subset \rad(J_{p^h, p^h}).\end{equation}
\noindent
Moreover, according to Lemma \ref{lemma3}, since char\,$R=p$, $x_0y_{p^h}-x_1y_{p^h-1}$, $x_{p^h-1}y_1-x_{p^h}y_0\in\rad(J_{p^h})$, so that
\begin{equation}\label{5a} J_{p^h, p^h}\subset\rad(J_{p^h}). \end{equation}
\noindent
The claim follows from (\ref{4a}), (\ref{3a}) and (\ref{5a}). \par\smallskip\noindent
We have just proven that in characteristic $p$, ideal $I_{p^h, p^h}$ can be generated, up to radical, by $2{p^h\choose 2}+1=p^h(p^h-1)+1$ minors of $A_{p^h, p^h}$.
 
\section{Set-theoretic complete intersections}
Now suppose that $K$ is an algebraically closed field, and that the set \underline{$x$}$=\{x_{ij}\}$ of entries of matrix $A$ defined in (\ref{2})  is a set of $N=c+d+2$ indeterminates over $K$ (here $x_{ij}\ne x_{hk}$ for $(i,j)\ne(h,k)$). Let $I$ be the ideal of $R=K[\mbox{\underline{$x$}}]$ generated by the 2-minors of $A$. Then $I$ is a prime ideal of height $c+d-1$.  The variety $V(I)\subset {\bf P}^{N-1}$ it defines is called a {\it rational normal scroll}.  See \cite{H} for the basic notions. \newline
Suppose that $P_1,\dots, P_s$ are elements of $R$ generating $I$ up to radical. It is well-known that $s\geq\,$height($I$). Whenever we can take $s=\,$height($I$) (or, more generally, $s\leq\,$height($I$)\,$+1$), ideal $I$ is called a {\it set-theoretic complete intersection} (or an {\it almost set-theoretic complete intersection}) on $P_1,\dots, P_s$.  
 
\begin{example}\label{example1} {\rm According to Corollary \ref{corollary4}, if char\,$K=2$, the ideal $I_{2,2}$ of $R=K[x_0, x_1, x_2, y_0, y_1, y_2]$ generated by the 2-minors of matrix
$$A=\left(\begin{array}{ccccc} x_0&x_1&\vline& y_0&y_1\\
x_1&x_2&\vline& y_1&y_2
\end{array}\right)$$
\noindent is generated, up to radical, by the following three binomials:
$$x_0x_2-x_1^2,\ y_0y_2-y_1^2,\ x_0y_2-x_2y_0.$$
\noindent
Since ht\,$I_{2,2}=3$, it follows that $I_{2,2}$ is a set-theoretic complete intersection on these three binomials.}
\end{example}
We will show that the property presented in Example \ref{example1} is shared by a larger class of rational normal scrolls. We first quote a result by Verdi \cite{V} (see also \cite{RV}, Section 1). For all $i=1,\dots, c-1$, let
$$F_i=\sum_{k=0}^i(-1)^k{i\choose k}x_{i+1}^{i-k}x_kx_i^k.$$
\noindent
\begin{proposition}\label{proposition5a} For all $c\geq 1$
$$(S_c)=\rad(F_1,\dots, F_{c-1}).$$
\noindent
In particular, ideal $(S_c)$ is a set-theoretic complete intersection. 
\end{proposition}
Now, for all $i=1,\dots, d-1$, let
$$G_i=\sum_{k=0}^i(-1)^k{i\choose k}y_{i+1}^{i-k}y_ky_i^k.$$
\noindent
Then, by Proposition \ref{proposition5a}, $(T_d)=\rad(G_1,\dots, G_{d-1}).$ In view of Corollary \ref{corollary4}, 
we deduce the following
\begin{corollary}\label{corollary4a} Let $p$ be a prime and suppose that $K$ is a field of characteristic $p$. Then, for all positive integers $h$, 
$$I_{p^h, p^h}=\rad(F_1,\dots, F_{p^h-1}, G_1, \dots, G_{p^h-1}, x_0y_{p^h}-x_{p^h}y_0).$$
\noindent
In particular, $I_{p^h, p^h}$ is a set-theoretic complete intersection. 
\end{corollary}
We also recall a result, due to Bardelli and Verdi \cite{BV} (see \cite{RV}, Section 2), on the ideal of $2$-minors of matrix $A_{c,1}$. Set
$$F_c=\sum_{k=0}^c(-1)^k{c\choose k}y_1^{c-k}x_ky_0^k.$$
\noindent
  \begin{proposition}\label{proposition5} For all $c\geq2$, 
$$I_{c,1}=\rad(F_1,\dots, F_{c-1}, F_c).$$
\noindent
In particular, $I_{c,1}$ is a set-theoretic complete intersection.
\end{proposition}
As a consequence of the results of Section 2, Proposition \ref{proposition5} can be extended in the following way: 
\begin{corollary}\label{corollary6} Suppose that $c,d\geq 2$. Then
$$I_{c,d}=\rad(F_1,\dots, F_{c-1}, F_c, G_1, \dots, G_{d-1}, x_0y_d-x_1y_{d-1}).$$
\noindent
In particular, $I_{c,d}$ is an almost set-theoretic complete intersection.
\end{corollary}
\demo The inclusion $\supset$ is clear. On the other hand, by Theorem \ref{main}, 
\begin{eqnarray*} I_{c,d}&=&\rad((S_c)+ (T_d)+(x_{c-1}y_1-x_cy_0)+(x_0y_d-x_1y_{d-1}))\\
&\subset&\rad(I_{c,1}+(T_d)+(x_0y_d-x_1y_{d-1}))\\
&=&\rad(F_1, \dots, F_{c-1}, F_c, G_1, \dots, G_{d-1}, x_0y_d-x_1y_{d-1}),
\end{eqnarray*}
\noindent
where the last equality follows from Proposition \ref{proposition5}. This completes the proof. 
\begin{example}{\rm Consider the matrix 
$$A_{4,3}=\left(\begin{array}{cccccccc}
x_0 & x_1 & x_2 & x_3 &\vline& y_0 & y_1 & y_2\\
x_1 & x_2 & x_3 & x_4 &\vline& y_1 & y_2 & y_3
\end{array}\right).$$
\noindent
According to Corollary \ref{corollary6}, the ideal $I_{4,3}$ of $K[x_0, x_1, x_2, x_3, x_4, y_0, y_1, y_2, y_3]$ generated by the 2-minors of $A_{3,4}$ is generated, up to radical, by the following 7 polynomials:
\begin{eqnarray*}
F_1&=&x_0x_2-x_1^2\\
F_2&=&x_0x_3^2-2x_1x_2x_3+x_2^3\\
F_3&=&x_0x_4^3-3x_1x_3x_4^2+3x_2x_3^2x_4-x_3^4\\
F_4&=&x_0y_1^4-4x_1y_0y_1^3+6x_2y_0^2y_1^2-4x_3y_0^3y_1+x_4y_0^4\\
G_1&=&y_0y_2-y_1^2\\
G_2&=&y_0y_3^2-2y_1y_2y_3+y_2^3\\
H&=&x_0y_3-x_1y_2.
\end{eqnarray*}
}\end{example}

\end{document}